\documentclass[12pt,leqno]{article}

=155mm\textheight =220mm

\usepackage{amsmath}
\usepackage{amscd}
\usepackage{amsmath}
\usepackage{amssymb}
\usepackage{latexsym}
\makeatletter

\@addtoreset{equation}{section}
\makeatother

\newtheorem{them}{Theorem}
\newtheorem{thm}{Theorem}[section]
\newtheorem{pr}[thm]{Proposition}
\newtheorem{df}[thm]{Definition}
\newtheorem{lm}[thm]{Lemma}
\newtheorem{cor}[thm]{Corollary}

\newcommand{\qed}{\hfill{\rule{7pt}{7pt}}}

\input xy
\xyoption{all}
\CompileMatrices

\begin{document}
\title{Ramification of local fields\\
with imperfect residue fields III}
\author{{\sc Takeshi Saito}}
\maketitle

\begin{abstract}
The graded quotients
of the logarithmic 
ramification groups
of a local field
of mixed characteristic
is killed by the residue characteristic.
Its characters are
described by differential forms.
\end{abstract}

Let $K$ be a complete
discrete valuation field
and
$F$ be the residue field.
We do not assume that
$F$ is perfect.
We fix a separable closure
$\bar K$ of $K$.
The residue field $\bar F$
of $\bar K$ is an algebraic closure
of $F$.
In \cite[Definition 3.12]{AS1},
we defined a filtration
by ramification groups
on the absolute Galois group
$G_K={\rm Gal}(\bar K/K)$
and a logarithmic variant.
In this paper,
we only consider the
logarithmic variant
and, by dropping the suffix log
in the notation, let
$(G_K^r)_{r\in {\mathbb Q},r>0}$
denote the decreasing filtration
by logarithmic ramification
groups.
For $r>0$,
we put
$G_K^{r+}=
\overline{\bigcup_{s>r}G_K^s}$.
In \cite[Theorem 5.12]{AS2},
we proved that
the graded quotient
${\rm Gr}^rG_K=
G_K^r/G_K^{r+}$
is an abelian group.
In this paper, we prove the following.

\begin{them}\label{them1}
Assume that the residue field
$F$ is of characteristic $p>0$.
Then,
the abelian group 
${\rm Gr}^rG_K=
G_K^r/G_K^{r+}$
is annihilated by $p$
for every $r>0$.
\end{them}

Liang Xiao claims Theorem \ref{them1}
in \cite[Theorem 3.7.3]{Xiao}.

Theorem \ref{them1} is proved 
in the equal characteristic case
in \cite[Corollary 1.27]{Jus}.
Similarly as
in the equal characteristic case,
it is reduced to the case
where the residue field $F$
is of finite type over
a perfect subfield $k$.
In this case,
an $F$-vector space
$\Omega_F(\log)$ of finite dimension
fitting in an exact sequence
$0\to \Omega^1_{F/k}
\to \Omega_F(\log)
\to F\to 0$
is defined (see (\ref{eqOmF})).
Let ${\rm ord}_{\bar K}$
be the valuation
of $\bar K$
extending the normalized
valuation ${\rm ord}_K$ of $K$
and we put 
${\mathfrak m}^r_{\bar K}=
\{x\in \bar K\mid
{\rm ord}_{\bar K}x\ge r\}$
and
${\mathfrak m}^{r+}_{\bar K}
=
\{x\in \bar K\mid
{\rm ord}_{\bar K}x>r\}$
and we consider the
$\bar F$-vector space
$\Theta_{\bar F}^{(r)}=
{\rm Hom}_F(
\Omega_F(\log),
{\mathfrak m}^r_{\bar K}/
{\mathfrak m}^{r+}_{\bar K})$
as a smooth additive algebraic group
over $\bar F$.
In \cite[(5.12.1)]{AS2},
further
a canonical surjection
\begin{equation}
\pi_1(\Theta_{\bar F}^{(r)})^{\rm ab}
\to
{\rm Gr}^rG_K
\label{eqcan}
\end{equation}
is defined (see (\ref{eqcan1})).
Let $\pi_1(\Theta_{\bar F}^{(r)})^{\rm alg}$
denote the quotient of
$\pi_1(\Theta_{\bar F}^{(r)})^{\rm ab}$
classifying \'etale isogenies.
Then, 
$\pi_1(\Theta_{\bar F}^{(r)})^{\rm alg}$
is a profinite abelian group killed by 
$p={\rm char} F>0$
and the character group
${\rm Hom}(\pi_1(\Theta_{\bar F}^{(r)})^{\rm alg},
{\mathbb F}_p)$
is canonically identified
with the dual space
${\rm Hom}_{\bar F}(
{\mathfrak m}^r_{\bar K}/
{\mathfrak m}^{r+}_{\bar K},
\Omega_F(\log)\otimes\bar F)$,
by pulling-back the
Artin-Schreier covering
${\mathbf G}_a
\to
{\mathbf G}_a:
t\mapsto t^p-t$
by linear form.
The main theorem of this paper
is the following.

\begin{them}\label{them2}
We assume that
the residue field
$F$ is finitely
generated over
a perfect subfield $k$
of characteristic $p>0$.
Then, the canonical surjection
{\rm (\ref{eqcan})}
factors through the quotient
$\pi_1(\Theta_{\bar F}^{(r)})^{\rm alg}$.
Consequently,
the abelian group
${\rm Gr}^rG_K$ is killed by $p$
and there exists a canonical injection
\begin{equation}
{\rm Hom}({\rm Gr}^rG_K,
{\mathbb F}_p)
\to
{\rm Hom}_{\bar F}(
{\mathfrak m}^r_{\bar K}/
{\mathfrak m}^{r+}_{\bar K},
\Omega_F(\log)\otimes\bar F).
\end{equation}
\end{them}

Since Theorem \ref{them2} is also proved 
in the equal characteristic case
in \cite[Theorem 1.24]{Jus},
we prove it in the mixed
characteristic case in this
paper.
The basic idea 
of the proof
is the same
as in the equal characteristic
case.
However, in the mixed
characteristic case, 
the projections
that played
a crucial role
in the proof
in the equal characteristic case
are not defined 
as maps of schemes.
They are defined only as an
infinitesimal deformation
of a morphism
of discrete valuation fields
in the sense of
Definition \ref{dfinf}.
We show that an infinitesimal
deformation induces
a functor of Galois categories
in Section 2
and that they satisfy
a transitivity property
in Section 3.
We prove Theorem \ref{them2}
in Section 4.
In Section 1,
we briefly recall fundamental
constructions
in the definition
of filtrations by
ramification groups.

The research is partly
supported by
Grants-in-aid for
Scientific Research B-18340002
and
A-22244001.

\section{Brief review of ramification
theory}

Let $K$ be a complete
discrete valuation field.
We assume that
$K$ is of characteristic $0$
and the residue field
$F$ is of characteristic $p>0$.
We assume
that $F$ is finitely
generated over
a perfect subfield $k$.
We define an $F$-vector
space of finite dimension
$\Omega_F(\log)$
by
\begin{equation}
\Omega_F(\log)
=
\left.
\bigl(\Omega^1_{F/k}
\oplus (F\otimes K^\times)\bigr)
\right/
(d \bar a- \bar
a\otimes a;
a\in {\cal O}_K^\times),
\label{eqOmF}
\end{equation}
by an abuse of notation
because $\Omega_F(\log)$
depends not only on $F$
but also on $K$.
It fits in an exact sequence
$0\to \Omega^1_{F/k}
\to \Omega_F(\log)
\to F\to 0$.
For $a\in K^\times$,
the image of
$1\otimes a$
is denoted by $d\log a$.

Let $\bar K$
be an algebraic closure
of $K$.
The residue field
$\bar F$ of $\bar K$
is an algebraic closure of $F$.
Let $G_K$ and $G_F$
be the absolute Galois
groups ${\rm Gal}(\bar K/K)$
and ${\rm Gal}(\bar F/F)$.
We have a canonical surjection
$G_K\to G_F$.
Let
$(G_K^r)_{r\in {\mathbb Q},r>0}$
denote the decreasing filtration
by logarithmic ramification
groups.
For $r>0$,
we put
$G_K^{r+}=
\overline{\bigcup_{s>r}G_K^s}$.
For a finite \'etale
$K$-algebra $L$,
we say that the log ramification
of $L$ is bounded by $r+$
if the natural action of
$G_K$ on the finite set
${\rm Hom}_K(L,\bar K)$
factors through the quotient
$G_K^{\le r+}=G_K/G_K^{r+}$.
Let 
${\cal C}_K^{\le r+}$
denote the category
of
finite \'etale
$K$-algebras
of log ramification
bounded by $r+$.
We identify
the category
${\cal C}_K^{\le r+}$
with the category
$(G_K^{\le r+}\text{-Sets})$
of finite sets with
continuous action of
$G_K^{\le r+}$
by the natural anti-equivalence
defined by
the fiber functor
attaching 
${\rm Hom}_K(L,\bar K)$
to $L$.

In the following of this section,
we assume that
$r>0$ is an integer. 
Let $\Theta^{(r)}$
denote the $F$-vector
space
${\rm Hom}_F(\Omega_F(\log),
{\mathfrak m}_K^r/
{\mathfrak m}_K^{r+1})$
of finite dimension
regarded as a smooth
algebraic group over $F$.
We consider a natural action
of $G_F$ on
$\Theta^{(r)}_{\bar F}=
\Theta^{(r)}\times_F\bar F$.
Let 
$(G_K^{\le r+}\text{-}{\rm FEt}/\Theta^{(r)}_{\bar F})$
denote the category
of finite \'etale schemes
over $\Theta^{(r)}_{\bar F}$
with a continuous action
of $G_K^{\le r+}=
G_K/G_K^{r+}$
compatible with that of 
$G_F$ on $\Theta^{(r)}_{\bar F}$.
We briefly
recall the construction
of the functor
\begin{equation}
X_K^{(r)}\colon
{\cal C}_K^{\le r+}
\longrightarrow
(G_K^{\le r+}\text{-}{\rm FEt}/\Theta^{(r)}_{\bar F})
\label{eqXr}
\end{equation}
in \cite{AS2}
with a slight modification
replacing complete
local rings by schemes
of finite type.

\begin{lm}\label{lmP0}
Let $L=\prod_jL_j$ be a finite 
\'etale $K$-algebra
and we put
$S={\rm Spec}\ 
{\cal O}_K$ and
$T={\rm Spec}\ 
{\cal O}_L$.
Then, there exists
a commutative diagram
\begin{equation}
\begin{CD}
T@>{i'}>> Q_0@<<<E_0\\
@VVV @VVV@VVV\\
S@>i>> P_0@<<<D_0
\end{CD}
\label{eqP0}
\end{equation}
of schemes over
the ring $W(k)$
of Witt vectors
satisfying the following
conditions:
\begin{itemize}
\item[{\rm (\ref{eqP0}.1)}]
The schemes $P_0,Q_0,D_0$
and $E_0$ are smooth
over $W(k)$
and $D_0\subset P_0$
and $E_0\subset Q_0$
are divisors.
The vertical
arrows
are finite
and flat.
The left square is
cartesian.
\item[{\rm (\ref{eqP0}.2)}]
Let $s={\rm Spec}\ F\in S$
denote the closed point
and $t_j={\rm Spec}\ F_j
\in T$
denote the closed points.
Then, the maps $i$ and $i'$
induces isomorphisms
$\kappa(i(s))\to F$ and
$\kappa(i'(t_j))\to F_j$
of residue fields.
The closed subschemes
$S\times_{P_0}D_0$
and 
$T\times_{Q_0}E_0$
are equal
to ${\rm Spec}\ F$
and 
to the reduced part
$(T\times_S
{\rm Spec}\ F)_{\rm red}
=
\coprod_j
{\rm Spec}\ F'_j$
respectively.
The canonical maps
$\Omega^1_{P_0/W(k)}(\log D_0)
\otimes F
\to \Omega_F(\log)$
and
$\Omega^1_{Q_0/W(k)}(\log E_0)
\otimes F'_j
\to \Omega_{F'_j}(\log)$
are isomorphisms
for every $j$.
On a neighborhood of $i'(t_j)$,
the pull-back
$D_0\times_{P_0}Q_0$
is equal to the divisor
$e_jE_0$ where $e_j$
is the ramification index,
for every $j$.
\end{itemize}
\end{lm}
{\it Proof.}
We take elements
$a_1,\ldots,a_n
\in {\cal O}_K$
such that the images
$\bar a_1,\ldots,
\bar a_n$ in $F$
form a transcendental basis
$F$ over $k$
and that
$F$ is a finite separable
extension of 
$k(\bar a_1,\ldots,\bar a_n)$.
Let $A_0$
be the henselization of
the subring
$W(k)[a_1,\ldots,a_n]
\subset {\cal O}_K$
at the prime ideal $(p)$
and $K_0$ be the fraction 
field of the completion of
the henselian
discrete valuation ring $A_0$.

Then $K$ is a finite separable
extension of $K_0$.
Let $K_1\subset K$
be the maximum unramified 
subextension over $K_0$.
Then, there exist
unique finite flat
normal $A_0$-subalgebras
$A_1\subset A$ of ${\cal O}_K$ 
such that
the natural maps
$A\otimes_{A_0}
{\cal O}_{K_0}
\to {\cal O}_K$ and
$A_1\otimes_{A_0}
{\cal O}_{K_0}
\to {\cal O}_{K_1}$ 
are isomorphisms.

We take a prime element $\pi$ of
$A$ and let $f\in A_1[t]$
be the minimal polynomial.
Let $A_1\{t\}$
be the henselization
at the maximal ideal $(p,t)$.
Then, we obtain an isomorphism
$A_1\{t\}/(f)\to A$.
It induces an isomorphism
$A_1\{t\}/(f,t)\to F$.

We may assume $L$
is a finite separable
extension of $K$.
Similarly, there exists
a unique finite flat
normal $A_0$-subalgebra
$B$ of ${\cal O}_L$ 
such that
the natural map
$B\otimes_{A_0}
{\cal O}_{K_0}
\to {\cal O}_L$ is
an isomorphism.

Let $F'$ be the
residue field of $L$.
We take elements
$b_1,\ldots,b_n
\in B$
such that the images
$\bar b_1,\ldots,
\bar b_n$ in $F'$
form a transcendental basis
$F'$ over $k$
and that
$F'$ is a finite separable
extension of 
$k(\bar b_1,\ldots,\bar b_n)$.
Let $B_0$
be the henselization of
the subring
$W(k)[b_1,\ldots,b_n]
\subset {\cal O}_L$
at the prime ideal $(p)$.
Then, we obtain
$L_0\subset L_1\subset L$
and $B_0\subset B_1\subset B$
as above.
We take a prime element $\pi'$ of
$B$ and let $g\in B_1[t']$
be the minimal polynomial.
Then, we obtain an isomorphism
$B_1\{t'\}/(g)\to B$.
It induces an isomorphism
$B_1\{t'\}/(g,t')\to F'$.

Since $A_1$ is essentially smooth
over $W(k)$,
there exists a map
$A_1\to B_1\{t'\}$
over $W(k)$
lifting the composition
$A_1\to A\to B=B_1\{t'\}/(g)$.
We put $\pi=u\pi^{\prime e}$
for $u\in B^\times$
and take a lifting
$\tilde u\in B_1\{t'\}^\times$.
We extend the map
$A_1\to B_1\{t'\}$
to a map
$A_1\{t\}\to B_1\{t'\}$
by sending $t$ to
$\tilde u\cdot t^{\prime e}$.
Thus, we obtain a commutative
diagram
\begin{equation}
\begin{CD}
B@<<< B_1\{t'\}@>{t'\mapsto 0}>> B_1\\
@AAA @AAA @AAA\\
A@<<< A_1\{t\}@>{t\mapsto 0}>> A_1
\end{CD}
\label{eqAB}
\end{equation}
of $W(k)$-algebras.

We show that the left square
gives an isomorphism
$A\otimes_{A_1\{t\}}
B_1\{t'\}\to B$.
Since the maximal
ideal of $A_1$
is generated by the 
image of $f$,
the maximal
ideal of $B_1$
is also generated by the 
image of $f$.
Hence, the image of
$f$ in $B_1\{t'\}$
is not in the square
of the maximal ideal
and we have
$(f)=(g)$ as ideals of
$B_1\{t'\}$.
Therefore the map
$A\otimes_{A_1\{t\}}
B_1\{t'\}
=B_1\{t'\}/(f)
\to B=B_1\{t'\}/(g)$
is an isomorphism.
Consequently,
the map $A_1\{t\}
\to B_1\{t'\}$
is finite flat.
Since the question
is \'etale local on
neighborhoods
of the images of
$S$ and $T$,
we deduce a diagram
(\ref{eqP0})
satisfying the conditions
(\ref{eqP0}.1) and 
(\ref{eqP0}.2)
from
the diagram (\ref{eqAB}).
\qed

We define a modification
$P^{(r)}$
of the scheme
$P_0\times_{W(k)}S$
as follows.
We take a blow-up
of $P_0\times_{W(k)}S$
at $D_0\times_{W(k)}{\rm Spec}\ F$
and define 
a scheme $P$ over $S$
to be the complement
of the union
of the proper transforms
of $P_0\times_{W(k)}{\rm Spec}\ F$
and 
$D_0\times_{W(k)}S$.
The map $S\to P_0$
induces a section
$S\to P$.
We regard
$S_r={\rm Spec}\
{\cal O}_K/
{\mathfrak m}^r_K$
as a closed subscheme
of $P$
by the composition
$S_r\to S\to P$.
We consider
the blow-up of $P$
at the closed subscheme
$S_r$
and define
$P^{(r)}$
to be the complement
of the proper
transform
of the closed
fiber $P\times_S
{\rm Spec}\ F$.
The schemes $P$
and $P^{(r)}$
are smooth over $S$.

More concretely,
the schemes $P$
and $P^{(r)}$
are described as follows.
Assume $P_0={\rm Spec}\ A_0$
is affine and the divisor
$D_0$ is defined
by $t\in A_0$.
The image $\pi
\in {\cal O}_K$
of $t$
by the map
$A_0\to {\cal O}_K$
corresponding to 
$S\to P_0$
is a uniformizer of $K$.
Then, we have
$P={\rm Spec}\ A$
for $A=A_0\otimes_{W(k)}
{\cal O}_K[U^{\pm 1}]
/(Ut-\pi)$.
Let $I$ be the kernel
of the surjection
$A\to {\cal O}_K$
induced by
$A_0\to {\cal O}_K$
and $U\mapsto 1$.
Then, we have
$P^{(r)}={\rm Spec}\ A^{(r)}$
for $A^{(r)}=
A[I/\pi^r]
\subset A[1/\pi]$.

The closed fiber
$P^{(r)}_F=
P^{(r)}\times_S{\rm Spec}\ F$
is canonically identified
with the affine space
$\Theta^{(r)}$
as follows.
The canonical map
$\Omega^1_{P_0/S}
\otimes_{{\cal O}_{P_0}}
{\cal O}_P
\to 
\Omega^1_{P/S}$
is uniquely extended
to an isomorphism
$\Omega^1_{P_0/S}(\log D_0)
\otimes_{{\cal O}_{P_0}}
{\cal O}_P
\to 
\Omega^1_{P/S}$.
Let ${\cal I}
\subset
{\cal O}_P$
denote
the ideal sheaf
defining
the closed subscheme $S\subset P$.
Then,
the closed fiber
$P^{(r)}_F$
is canonically
identified with
the $F$-vector space
${\rm Hom}_F
({\cal I}/{\cal I}^2
\otimes F,
{\mathfrak m}^r_K/
{\mathfrak m}^{r+1}_K)$
regarded as
an affine space.
By the isomorphisms
$\Omega^1_{P_0/S}(\log D_0)
\otimes_{{\cal O}_{P_0}}
{\cal O}_P
\to 
\Omega^1_{P/S},\
{\cal I}/{\cal I}^2
\to
\Omega^1_{P/S}
\otimes_{{\cal O}_P}
{\cal O}_S$
and 
$\Omega^1_{P_0/S}(\log D_0)
\otimes F
\to
\Omega_F(\log)$,
we obtain a canonical
isomorphism
\begin{equation}
P^{(r)}_F
\to
\Theta^{(r)}.
\label{eqth}
\end{equation}

Let $Q^{(r)}_{\bar S}$
be the normalization
of the base change
$Q_0\times_{P_0}
P^{(r)}_{\bar S}$
and
$Q^{(r)}_{\bar F}$
be the closed fiber.
Then, by the description
of log ramification groups
\cite[Section 5.1]{AS2},
the log ramification
of $L$ is bounded by $r+$
if and only if
the finite map
$Q^{(r)}_{\bar F}
\to 
\Theta^{(r)}_{\bar F}$
is \'etale.
Further, it is
shown in \cite[Lemma 5.10]{AS2} that,
if
the log ramification
of $L$ is bounded by $r+$,
the finite \'etale
scheme
$Q^{(r)}_{\bar F}
\to 
\Theta^{(r)}_{\bar F}$
with the natural action
of $G_K$ is
independent of the choice
of a diagram (\ref{eqP0})
and is well-defined
up to unique isomorphism.
The functor
$X_K^{(r)}\colon
{\cal C}_K^{\le r+}
\to
(G_K^{\le r+}\text{-}{\rm FEt}/\Theta^{(r)}_{\bar F})$
(\ref{eqXr})
is defined by attaching
$Q^{(r)}_{\bar F}$ to $L$.
The composition
of
$X_K^{(r)}\colon
{\cal C}_K^{\le r+}
\to
(G_K^{\le r+}\text{-}{\rm FEt}/\Theta^{(r)}_{\bar F})$
with the
fiber functor 
$F_{\bar 0}\colon
(G_K^{\le r+}\text{-}{\rm FEt}/\Theta^{(r)}_{\bar F})
\to
(G_K^{\le r+}\text{-Sets})$
at the
origin $0\in \Theta^{(r)}_{\bar F}$
recovers
the natural equivalence
of categories
${\cal C}_K^{\le r+}
\to
(G_K^{\le r+}\text{-Sets})$.

Further,
it is shown in
\cite[Theorem 5.12]{AS2} that, 
for a finite \'etale $K$-algebra
$L$ of log ramification bounded by
$r+$,
the finite \'etale covering
$X_K^{(r)}(L)\to \Theta^{(r)}_{\bar F}$
is trivialized by
a universal abelian covering
$\Theta^{(r) {\rm ab}}_{\bar F}$.
Thus, forgetting the Galois action
on $X_K^{(r)}(L)$ and taking
the fiber functor at the origin
$0\in \Theta^{(r)}_{\bar F}$,
we obtain a group homomorphism
$\pi_1(\Theta^{(r)}_{\bar F})^{\rm ab}
\to G_K^{\le r+}$
defined by the functor $X_K^{(r)}$.
It induces 
a canonical surjection
\begin{equation}
\pi_1(\Theta_{\bar F}^{(r)})^{\rm ab}
\to
{\rm Gr}^rG_K
\label{eqcan1}
\end{equation}
\cite[(5.12.1)]{AS2}.
It is compatible 
with the $G_K$-action
defined by the
actions of $G_F$
on $\Theta_{\bar F}^{(r)}$
and the conjugate action on
${\rm Gr}^rG_K$.
Since the inertia group
$I={\rm Ker}(G_K\to G_F)$
acts trivially
on $\Theta^{(r)}_{\bar F}$,
it follows
that ${\rm Gr}^rG_K$
is a central subgroup
of $I/G_K^{r+}$.

\section{Infinitesimal deformation}

Let $K$ and $K'$ be
complete discrete valuation fields.
We say that a morphism
$f\colon K\to K'$ of fields
is an extension of complete
discrete valuation fields
if it induces a flat local morphism
${\cal O}_K\to
{\cal O}_{K'}$,
also denoted by $f$
by abuse of notation.
The integer $e>0$
characterized by
$f({\mathfrak m}_K){\cal O}_{K'}
=
{\mathfrak m}^e_{K'}$
is called the ramification
index of $f$.

\begin{df}\label{dfinf}
Let $f\colon K\to K'$
be an extension
of complete discrete 
valuation fields
of ramification 
index $e$.
For an integer $r>0$,
we call the pair
$(f,\varepsilon)$
with an $f$-linear morphism
$\varepsilon\colon
\Omega_F(\log)\to
{\mathfrak m}_{K'}^{er}/
{\mathfrak m}_{K'}^{er+1}$
an infinitesimal
deformation of $f$.
\end{df}

We define a composition 
of infinitesimal deformations.
Let $f\colon K\to K',
g\colon K'\to K''$
be extensions
of complete discrete valuation fields
of ramification indices
$e,e'$
and $\varepsilon\colon
\Omega_F(\log)\to
{\mathfrak m}_{K'}^{er}/
{\mathfrak m}_{K'}^{er+1}$
and $\varepsilon'\colon
\Omega_{F'}(\log)\to
{\mathfrak m}_{K''}^{ee'r}/
{\mathfrak m}_{K''}^{ee'r+1}$
be $f$-linear and
$g$-linear morphisms.
Let $g_*\colon
{\mathfrak m}_{K'}^{er}/
{\mathfrak m}_{K'}^{er+1}\to
{\mathfrak m}_{K''}^{ee'r}/
{\mathfrak m}_{K''}^{ee'r+1}$
be the map induced by $g$
and
$f_*\colon
\Omega_F(\log)\to
\Omega_{F'}(\log)$
be the map induced by $f$.
Then, we define
the composition
$(g,\varepsilon')\circ
(f,\varepsilon)$
to be
$(g\circ f,
g_*\circ \varepsilon+
\varepsilon'\circ f_*)$.

Let $f\colon K\to K'$
be an extension
of complete discrete valuation
field of ramification
index $e$.
Let $r>0$ an integer
and we put $r'=er$.
Let $\varepsilon
\colon
\Omega_{{\cal O}_K}(\log)
\to {\mathfrak m}_{K'}^{r'}/
{\mathfrak m}_{K'}^{r'+1}$
be an $f$-linear morphism.
For an infinitesimal
deformation
$(f,\varepsilon)$ of $f$,
we define a functor
\begin{equation}
f_{\varepsilon*}\colon 
{\cal C}_K^{\le r+}\to
{\cal C}_{K'}^{\le r'+}.
\label{eqf*}
\end{equation}

We take separable closures
$K\subset \bar K,
K'\subset \bar K'$
and an embedding 
$\bar f\colon \bar K
\to \bar K'$.
The residue fields
$\bar F, \bar F'$
of 
$\bar K, \bar K'$
are algebraic closures
of $F$ and of $F'$.
By \cite[Proposition 3.15 (3)]{AS1},
the map $f^*\colon
G_{K'}\to G_K$
induces
$G_{K'}^{\le er+}
\to G_K^{\le r+}$.
The embedding
$\bar f\colon \bar K
\to \bar K'$
induces
$\bar f\colon \bar F
\to \bar F'$
and defines a commutative
diagram
\begin{equation}
\begin{CD}
G_{K'}^{\le r'+}
@>{f^*}>> 
G_K^{\le r+}\\
@VVV @VVV\\
G_{F'}@>>> G_F
\end{CD}\label{eqGF}
\end{equation}
We recalled the
definition of a functor
$X_K^{(r)}\colon
{\cal C}_K^{\le r+}
\to 
(G_K^{\le r+}\text{-}{\rm FEt}/\Theta^{(r)}_{\bar F})$
(\ref{eqXr}) in Section 1.
The map $\varepsilon$
defines a geometric point
$\bar \varepsilon\colon
\bar F'\to \Theta^{(r)}_{\bar F}$.
The map
$\bar \varepsilon\colon
\bar F'\to \Theta^{(r)}_{\bar F}$
is compatible with
the morphism
$G_{F'}\to G_F$.
By the commutative diagram
(\ref{eqGF}),
the fiber functor
$F_{\bar \varepsilon}$
defines a functor
$(G_K^{\le r+}\text{-}{\rm FEt}/\Theta^{(r)}_{\bar F})
\to (G_{K'}^{\le r'+}\text{-Sets})\simeq
{\cal C}_{K'}^{\le r'+}$.
We define a functor
$f_{\varepsilon*}$
as the composition
\begin{equation}
F_{\bar \varepsilon}
\circ
X^{(r)}_K\colon
{\cal C}_K^{\le r+}\to
(G_K^{\le r+}\text{-}{\rm FEt}/\Theta^{(r)}_{\bar F})
\to {\cal C}_{K'}^{\le r'+}.
\label{eqfe}
\end{equation}

To describe a morphism
$f_{\varepsilon}^*\colon
G_{K'}^{\le r'+}\to
G_K^{\le r+}$
corresponding
to the functor
$f_{\varepsilon*}$,
we introduce a terminology.
Let $G$ and $G'$ be groups
and $C\subset G$
be a central subgroup.
For morphisms of groups
$\varphi\colon G'\to G$ and
$\psi\colon G'\to C$,
we call the morphism
$\varphi_\psi\colon
G'\to G$ 
defined by 
$\varphi_\psi(g)
=\varphi(g)\psi(g)$
for 
$g\in G'$,
the deformation
of
$\varphi\colon G'\to G$ by
$\psi\colon G'\to C$.

We consider the composition
\begin{equation}
\begin{CD}
G_{K'}^{\le r'+}
\to G_{F'}^{\rm ab}
@>{\varepsilon_*}>>
\pi_1(\Theta^{(r)}_{\bar F})^{\rm ab}
@>{\rm(\ref{eqcan1})}>>
{\rm Gr}^rG_K
\subset G_K^{\le r+}=
G_K/G_K^{r+}.
\end{CD}
\label{eqep*}
\end{equation}
As is remarked after
(\ref{eqcan1}),
the subgroup
${\rm Gr}^rG_K
\subset G_K^{\le r+}$
is a central subgroup
if the residue field
$F$ is separably closed.

\begin{lm}\label{lmhom}
Assume that the residue field
$F$ is separably closed.
Then, the functor
$f_{\varepsilon*}\colon 
{\cal C}_K^{\le r+}\to
{\cal C}_{K'}^{\le r'+}$
is compatible with
the deformation
$f^*_{\varepsilon^*}\colon
G_{K'}^{\le r'+}\to G_K^{\le r+}$
of 
$f^*\colon
G_{K'}^{\le r'+}\to G_K^{\le r+}$
by $\varepsilon_*\colon
G_{K'}^{\le r'+}\to {\rm Gr}^rG_K
\subset G_K^{\le r+}$
{\rm (\ref{eqep*})}.
\end{lm}

{\it Proof.}
We take a lifting
${\rm Spec}\ \bar F'
\to \Theta^{(r) {\rm ab}}_{\bar F}$
to a universal abelian
covering
of the geometric point
$\bar\varepsilon\colon
{\rm Spec}\ \bar F'
\to \Theta^{(r)}_{\bar F}$
and consider the bijection
\begin{equation}
X_K^{(r)}(L)
\times_{\Theta^{(r)}_{\bar F}}
{\rm Spec}\ \bar F'
\to \pi_0(X_K^{(r)}(L)
\times_{\Theta^{(r)}_{\bar F}}
\Theta^{(r) {\rm ab}}_{\bar F})
\label{eqbij}
\end{equation}
of finite sets
for a finite \'etale $K$-algebra
$L$ of ramification
bounded by $r+$.
By the definition of
the functor $f_{\varepsilon *}$,
the finite $G_{K'}^{\le r'+}$-set
$f_{\varepsilon *}(L)$
is defined as
$X_K^{(r)}(L)
\times_{\Theta^{(r)}_{\bar F}}
{\rm Spec}\ \bar F'$.
The bijection
(\ref{eqbij})
is compatible with
the map
$(f^*,\varepsilon_*)
\colon
G_{K'}^{\le r'+}
\to 
G_K^{\le r+}
\times
\pi_1(\Theta_{\bar F})^{\rm ab}$.
By the definition of
the canonical map
(\ref{eqcan1}),
the action of
$\pi_1(\Theta_{\bar F})^{\rm ab}$
on the finite set
$\pi_0(X_K^{(r)}(L)
\times_{\Theta^{(r)}_{\bar F}}
\Theta^{(r) {\rm ab}}_{\bar F})$
is the same as that
induced from the action
of $G_K^{\le r+}$
by (\ref{eqcan1}).
Thus the assertion 
follows.
\qed

If we choose a morphism of
fiber functors,
the functor
$f_{\varepsilon*}\colon 
{\cal C}_K^{\le r+}\to
{\cal C}_{K'}^{\le r'+}$ induces
a morphism of groups
$f_{\varepsilon}^*\colon
G_{K'}^{\le r'+}\to
G_K^{\le r+}$.
Without choosing
a morphism of
fiber functors,
it is still
well-defined up to conjugate.
Hence, for a representation $V$ of
$G_K^{\le r+}$
the restriction
${\rm Res}_{f,\varepsilon}V$
is well-defined up to
an isomorphism
as a representation of
$G_{K'}^{\le r'+}$.

\begin{cor}\label{corchi}
Assume that the residue field $F$
is separably closed.
Let $V$ be a represention of
$G_K^{\le r+}$ such that
the restriction to
${\rm Gr}^rG_K$ is a character
$\chi$.
Let $\varepsilon^*(\chi)$
denote the character
of $G_{K'}^{\le r'+}$
defined as the pull-back by
$\varepsilon_*
\colon G_{K'}^{\le r'+}
\to
{\rm Gr}^rG_K$
{\rm (\ref{eqep*})}.

Then, we have
an isomorphism
\begin{equation}
{\rm Res}_{f,\varepsilon}V
\to
{\rm Res}_{f}V
\otimes 
\varepsilon^*(\chi)
\end{equation}
of representations of
$G_{K'}^{\le r'+}$.
\end{cor}

\section{Transitivity}

Let $f\colon K\to K'$
be an extension of complete
discrete valuation fields.
We say that $K'$
is a smooth extension of $K$
if the ramification index is $1$
and if the residue field
$F'$ of $K'$ is a
finitely generated separable extension
of the residue field
$F$ of $K$.

Let $f\colon K\to K'$
be a smooth extension of 
complete discrete valuation
field and
let $\varepsilon
\colon \Omega_F(\log)
\to {\mathfrak m}^r_{K'}/
{\mathfrak m}^{r+1}_{K'}$
be an $f$-linear map.
We consider the dual
$${\rm Hom}_{F'}(
\Omega_{F'}(\log),
{\mathfrak m}^r_{K'}/
{\mathfrak m}^{r+1}_{K'})
\to
{\rm Hom}_F(
\Omega_F(\log),
{\mathfrak m}^r_K/
{\mathfrak m}^{r+1}_K)
\otimes_F{F'}$$
of the map
$\Omega_F(\log)
\otimes_FF'
\to
\Omega_{F'}(\log)$
induced by $f$.
We also consider
the translation
$+\varepsilon$
as a morphism
$\Theta^{\prime(r)}_{F'}=
{\rm Hom}_{F'}(
\Omega_{F'}(\log),
{\mathfrak m}^r_{K'}/
{\mathfrak m}^{r+1}_{K'})
\to \Theta^{\prime(r)}_{F'}$.
Their composition
defines a morphism of
schemes
$f_\varepsilon^*\colon
\Theta^{\prime(r)}_{\bar F'}
\to \Theta^{(r)}_{\bar F}$
compatible with
$G_{F'}\to G_F$
and hence
the pull-back functor
$f_{\varepsilon*}\colon
(G_K^{\le r+}\text{-}{\rm FEt}/\Theta^{(r)}_{\bar F})
\to
(G_{K'}^{\le r+}\text{-}
{\rm FEt}/\Theta^{\prime(r)}_
{\bar F'}).$

\begin{pr}\label{prsm}
Assume $f\colon K\to K'$
is smooth
and consider the diagram
$$\begin{CD}
{\cal C}_K^{\le r+}
@>{f_{\varepsilon *}}>>
{\cal C}_{K'}^{\le r+}\\
@V{X^{(r)}_K}VV 
@VV{X^{(r)}_{K'}}V\\
(G_K^{\le r+}\text{-}{\rm FEt}/\Theta^{(r)}_{\bar F})
@>{f_{\varepsilon *}}>>
(G_{K'}^{\le r+}\text{-}
{\rm FEt}/\Theta^{\prime(r)}_
{\bar F'})
\end{CD}$$
of functors.
Then, there exists
an isomorphism
\begin{equation}
f_{\varepsilon *}\circ
X^{(r)}_K
\to
X^{(r)}_{K'}
\circ
f_{\varepsilon *}
\label{eqsm}
\end{equation}
of functors.
\end{pr}

{\it Proof.}
We regard $\varepsilon\colon
\Omega_F(\log)
\to {\mathfrak m}_{K'}^r/
{\mathfrak m}_{K'}^{r+1}$
as an $F'$-rational point
$\varepsilon\colon
{\rm Spec}\ F'
\to \Theta^{(r)}_F$.
Let $L$ be a finite \'etale
algebra over $K$
of log ramification
bounded by $r+$.
We take a cartesian
diagram (\ref{eqP0})
over $W(k)$
as in Lemma \ref{lmP0}.
Take a morphism
$S'\to P^{(r)}$
over $S$
lifting
$\varepsilon\colon
{\rm Spec}\ F'
\to \Theta^{(r)}
\subset P^{(r)}$
and consider the
composition
$S'\to P^{(r)}
\to P\to P_0$.

We put $T'=S'\times_{P_0}Q_0$.
Since $\bar Q^{(r)}_{\bar F}
\to \bar P^{(r)}_{\bar F}=
\Theta^{(r)}_{\bar F}$
is \'etale,
the base chage
$Q_0\times_{P_0}P^{(r)}
\to P^{(r)}$
is \'etale on 
the complement of
the closed fiber in
a neighborhood of
the closed fiber.
Hence,
the $K'$-algebra
$L'=
\Gamma(T'\times_{S'}{\rm Spec}\ K',
{\cal O})$
is \'etale.
The fiber product
$T'\times_{Q_0}E_0$
is isomorphic to
$S'\times_{P_0}E_0
=
(S'\times_{P_0}D_0)
\times_{D_0}E_0
=
F'\times_F
(F\times_{D_0}E_0)
=
F'\times_F
(T\times_{Q_0}E_0)_{\rm red}$
and is reduced
since $F'$ is separable over $F$.
Hence,
$T'$ is the spectrum
of the integer ring ${\cal O}_{L'}$.

By the assumption
that $K'$ is smooth over $K$,
there exists a commutative diagram
$$\begin{CD}
S'@>>> P'_0@<<< D'_0\\
@VVV @VVV@VVV\\
S@>>> P_0@<<<D_0
\end{CD}$$
of schemes over $W(k)$
satisfying the following
conditions:
The vertical arrow
$P'_0\to P_0$ is smooth,
the right square is cartesian,
$S'\times_{P'_0}D'_0=
{\rm Spec}\ F'$
and $\Omega^1_{P'_0/W(k)}(\log D'_0)
\otimes F'\to \Omega_{F'}(\log)$
is an isomorphism.

We consider the diagram
\begin{equation}
\begin{CD}
T'@>>> Q'_0
@<<< E'_0\\
@VVV @VVV@VVV\\
S'@>>> P'_0@<<<D'_0
\end{CD}
\label{eqP'0}
\end{equation}
where the right square
is the base change
of that of (\ref{eqP0})
by $P'_0\to P_0$
and the left square
is cartesian.
Since 
$T'\times_{Q'_0}E'_0
=
T'\times_{Q_0}E_0
=
F'\times_F
(T\times_{Q_0}E_0)_{\rm red}
=
(T'\times_{Q_0}E_0)_{\rm red}$,
the diagram
(\ref{eqP'0}) satisfies the conditions
corresponding to
(\ref{eqP0}.1) and
(\ref{eqP0}.2).

We define
$Q^{(r)}_{\bar S}
\to P^{(r)}_{\bar S}$
and 
$Q^{\prime (r)}_{\bar S'}
\to P^{\prime(r)}_{\bar S'}$
and we identify
$P^{(r)}_{\bar S}=
\Theta^{(r)}_{\bar F}$
and
$P^{\prime (r)}_{\bar S'}=
\Theta^{(r)}_{\bar F'}$
as in (\ref{eqth}).
Then, the map
$P^{\prime (r)}_{\bar S'}
\to 
P^{(r)}_{\bar S}
\times_{\bar S}{\bar S'}$
induced by $P'\to P$
is smooth and hence the diagram
\begin{equation}
\begin{CD}
Q^{\prime (r)}_{\bar S'}
@>>>
Q^{(r)}_{\bar S}
\\
@VVV@VVV\\
P^{\prime (r)}_{\bar S'}
@>>>
P^{(r)}_{\bar S}
\end{CD}
\label{eqPP'}
\end{equation}
is cartesian.
Since the diagram
\begin{equation}
\begin{CD}
\Theta^{\prime (r)}_{\bar F'}
@>{f_{\varepsilon*}}>>
\Theta^{(r)}_{\bar F}
\\
@VVV@VVV\\
P^{\prime (r)}_{\bar S'}
@>>>
P^{(r)}_{\bar S}
\end{CD}
\label{eqPT'}
\end{equation}
is cartesian,
we obtain a cartesian diagram
\begin{equation}
\begin{CD}
X_{K'}^{(r)}(L')
@>>> X_K^{(r)}(L)
\\
@VVV@VVV\\
\Theta^{\prime (r)}_{\bar F'}
@>{f_{\varepsilon*}}>>
\Theta^{(r)}_{\bar F}
\end{CD}
\label{eqXT'}
\end{equation}
compatible with the
group homomorphism
$G_{K'}^{\le r+}\to G_K^{\le r+}$.

Thus by the definition
of the functor
$f_{\varepsilon *}$, we obtain
an isomorphism
$L'\to f_{\varepsilon *}(L)$.
The diagram
(\ref{eqXT'})
defines an isomorphism
$X_{K'}^{(r)}(L')
\to f_{\varepsilon*} X_K^{(r)}(L)$
in the category
$(G_{K'}^{\le r+}\text{-}{\rm FEt}/
\Theta^{\prime (r)}_{\bar F'})$.
The isomorphism is functorial in $L$
and they define an isomorphism
$X_{K'}^{(r)}\circ f_{\varepsilon*}
\to f_{\varepsilon*} 
\circ X_K^{(r)}$
of functors.
\qed

We deduce the following
transitivity.

\begin{cor}\label{cortra}
Let $f\colon K\to K'$
and $g\colon K' \to K''$
be extensions
of complete discrete valuation
fields.
We assume $f$ is {\rm smooth}
and let $e'$ be the ramification
index of $K''$ over $K'$.
Let $\varepsilon\colon
\Omega_F(\log)\to
{\mathfrak m}_{K'}^r/
{\mathfrak m}_{K'}^{r+1}$
and $\varepsilon'\colon
\Omega_{F'}(\log)\to
{\mathfrak m}_{K''}^{e'r}/
{\mathfrak m}_{K''}^{e'r+1}$
be infinitesimal deformations
and we put
$(g,\varepsilon')\circ
(f,\varepsilon)
=(g\circ f,\varepsilon'')$.
Then,
there exists an isomorphism
of functors:
\begin{equation}
(g\circ f)_{
\bar \varepsilon''*}
\to
g_{\varepsilon'*}
\circ
f_{\varepsilon*}.
\label{eqtra}
\end{equation}
\end{cor}

{\it Proof.}
The composition 
of the morphism
(\ref{eqsm})
with the fiber functor
$F_{\bar \varepsilon'}$
gives a morphism
$$F_{\bar \varepsilon'}
\circ
f_{\varepsilon *}\circ
X^{(r)}_K
\to
F_{\bar \varepsilon'}
\circ
X^{(r)}_{K'}
\circ
f_{\varepsilon *}=
g_{\varepsilon' *}
\circ
f_{\varepsilon *}$$
of functors.
By the canonical 
isomorphism
$F_{\bar \varepsilon'}
\circ
f_{\varepsilon *}
\to
F_{\bar \varepsilon''}$,
the first term
$F_{\bar \varepsilon'}
\circ
f_{\varepsilon *}\circ
X^{(r)}_K$
is identified with
the functor
$(g\circ f)_{
\bar \varepsilon''*}$.
\qed

\section{Proof of Theorem \ref{them2}}

We prove Theorem \ref{them2}
in the introduction.
It is reduced to
the case where
$r>0$ is an integer,
by considering
the base change
by log smooth extension
as in \cite[Lemma 1.22]{Jus}.
We regard an
$F$-vector space $V$
of finite dimension
as a smooth additive
algebraic group over
$F$ and 
let ${\rm pr}_1,
{\rm pr}_2\colon
V\times V\to V$
be the projections
and $-\colon
V\times V\to V$
be the subtraction
$(x,y)\mapsto y-x$.
By Lemma \ref{lmalg}
below,
it suffices to prove
the following.

\begin{lm}\label{lmchi}
Let $\chi$
be a character of
${\rm Gr}^rG_K$
and regard it also
as a character
of $\pi_1(\Theta^{(r)}_{\bar F})^{\rm ab}$
by the surjection
{\rm (\ref{eqcan1})}.
Then, 
we have an equality
${\rm pr}_2^*\chi
=
{\rm pr}_1^*\chi
\cdot
-^*\chi$
of characters of
$\pi_1(\Theta^{(r)}_{\bar F}
\times \Theta^{(r)}_{\bar F}
)^{\rm ab}$.
\end{lm}

\begin{lm}[{\cite[Lemma 1.23]{Jus}}]\label{lmalg}
Let $F$ be an algebraically
closed field of 
characteristic $p>0$
and regard an
$F$-vector space $V$
of finite dimension
as a smooth additive
algebraic group over $F$.
Let $\pi_1(V)^{\rm alg}$
be the quotient
of the abelian fundamental
group $\pi_1(V)^{\rm ab}$
classifying
\'etale isogenies.
Then, for a character
$\chi$ of $\pi_1(V)^{\rm ab}$,
the following
conditions are equivalent:

{\rm (1)}
$\chi$ factors through
the quotient
$\pi_1(V)^{\rm alg}$.

{\rm (2)}
We have an equality
${\rm pr}_2^*\chi
=
{\rm pr}_1^*\chi
\cdot
-^*\chi$
of characters of
$\pi_1(V\times V)^{\rm ab}$.
\end{lm}

To prove Lemma \ref{lmchi},
we use the geometric construction
in Section 1.
We consider the 
smooth scheme 
$P^{(r)}$ over $S$
and the fiber product
$P^{(r)}\times_S
P^{(r)}$.
The closed fibers
$P^{(r)}_F$ 
and
$P^{(r)}_F\times_F
P^{(r)}_F$
are identified with
$\Theta^{(r)}_F$
and
$\Theta^{(r)}_F
\times \Theta^{(r)}_F$.
Let $\xi
\in 
\Theta^{(r)}_F
\subset
P^{(r)}$
and $\eta
\in 
\Theta^{(r)}_F
\times_F
\Theta^{(r)}_F
\subset
P^{(r)}\times_S
P^{(r)}$
be the generic points.
Define complete discrete
valuation fields
$K'$ and $K''$
to be the fraction fields
of the completions of
the local rings
${\cal O}_{P^{(r)},\xi}$
and 
${\cal O}_{P^{(r)}\times_S
P^{(r)},\eta}$
respectively.
They are smooth extensions
of $K$.
Let $1\colon K\to K'$
denote the canonical map
and ${\rm p}_1,{\rm p}_2
\colon K'\to K''$
denote the map
induced by the projections
$P^{(r)}\times_S
P^{(r)}\to
P^{(r)}$.

We define infinitesimal deformations
$\varepsilon\colon
\Omega_F(\log)\to 
{\mathfrak m}^r_{K'}/
{\mathfrak m}^{r+1}_{K'}$
and
$\varepsilon'\colon
\Omega_{F'}(\log)\to 
{\mathfrak m}^r_{K''}/
{\mathfrak m}^{r+1}_{K''}$.
The residue field $F'$
of $K'$
is the function field
of $\Theta^{(r)}_F$
and hence is the
fraction field
of the symmetric algebra
$S^\bullet_F
{\rm Hom}_F(
{\mathfrak m}^r_K/
{\mathfrak m}^{r+1}_K,
\Omega_F(\log))$
of the dual vector space.
We define
$\varepsilon\colon
\Omega_F(\log)\to 
{\mathfrak m}^r_{K'}/
{\mathfrak m}^{r+1}_{K'}
=F'\otimes_F
{\mathfrak m}^r_K/
{\mathfrak m}^{r+1}_K$
to be the map
\begin{align*}
\Omega_F(\log)\to &
{\rm Hom}_F(
{\mathfrak m}^r_K/
{\mathfrak m}^{r+1}_K,
\Omega_F(\log))
\otimes_F
{\mathfrak m}^r_K/
{\mathfrak m}^{r+1}_K\\
&\subset
S^\bullet_F
({\rm Hom}_F(
{\mathfrak m}^r_K/
{\mathfrak m}^{r+1}_K,
\Omega_F(\log)))
\otimes_F
{\mathfrak m}^r_K/
{\mathfrak m}^{r+1}_K
\subset
F'\otimes_F
{\mathfrak m}^r_K/
{\mathfrak m}^{r+1}_K
\end{align*}
where the arrow is
the inverse of the
isomorphism defined
by the evaluation.

Similarly,
the residue field $F''$
of $K''$
is the
fraction field
of the symmetric algebra
$S^\bullet_F
{\rm Hom}_F(
{\mathfrak m}^r_K/
{\mathfrak m}^{r+1}_K,
\Omega_F(\log)^{\oplus 2})$.
Since $K'$ is
a smooth extension of $K$,
we have an exact sequence
$0\to
\Omega_F(\log)
\otimes_FF'
\to
\Omega_{F'}(\log)
\to 
\Omega_{F'/F}\to 0$.
Let 
$\varepsilon'\colon
\Omega_{F'}(\log)\to 
{\mathfrak m}^r_{K''}/
{\mathfrak m}^{r+1}_{K''}
=F''\otimes_F
{\mathfrak m}^r_K/
{\mathfrak m}^{r+1}_K$
be a ${\rm p}_1$-linear map
such that
the restriction to
$\Omega_F(\log)$
is 
\begin{align*}
\Omega_F(\log)\to 
&
\Omega_F(\log)^{\oplus 2}
\to 
{\rm Hom}_F(
{\mathfrak m}^r_K/
{\mathfrak m}^{r+1}_K,
\Omega_F(\log)^{\oplus 2})
\otimes_F
{\mathfrak m}^r_K/
{\mathfrak m}^{r+1}_K\\
&\subset
S^\bullet_F
({\rm Hom}_F(
{\mathfrak m}^r_K/
{\mathfrak m}^{r+1}_K,
\Omega_F(\log)^{\oplus 2}))
\otimes_F
{\mathfrak m}^r_K/
{\mathfrak m}^{r+1}_K
\subset
F''\otimes_F
{\mathfrak m}^r_K/
{\mathfrak m}^{r+1}_K
\end{align*}
where the first arrow is
the map
$x\mapsto (-x,x)$
and the second arrow
is the inverse of the
isomorphism defined
by the evaluation.

Finally, we define a morphism
$\mu\colon K'\to K''$
over $K$.
The subtraction map
$\Theta_F^{(r)}\times_F
\Theta_F^{(r)}\to
\Theta_F^{(r)}$
is dominant and
induces a morphism
$F'\to F''$ over $F$.
We define
$\mu\colon K'\to K''$
to be a morphism
over $K$
lifting
the map $F'\to F''$.

\begin{lm}\label{lmmu}
We have
\begin{align}
{\rm p}_2\circ (1,\varepsilon)
&=
({\rm p}_1,\varepsilon')
\circ (1,\varepsilon)
\label{eq2e}
\\
\mu\circ (1,\varepsilon)
&=
({\rm p}_1,\varepsilon')
\circ 1.
\label{eqmue}
\end{align}
\end{lm}
{\it Proof.}
They are
equalities
of deformations
of the canonical map
$K\to K''$.
Hence, it suffices
to prove the equalities
of maps
$\Omega_F(\log)\to
{\mathfrak m}^r_{K''}
/{\mathfrak m}^{r+1}_{K''}
=F''\otimes
{\mathfrak m}^r_K
/{\mathfrak m}^{r+1}_K$.

For the left hand side of
(\ref{eq2e}),
it is induced by the map
$\Omega_F(\log)
\to
\Omega_F(\log)^{\oplus 2}$
sending $x$ to $(0,x)$.
For the right hand side of
(\ref{eq2e}),
it is induced by the sum
of the maps
$\Omega_F(\log)
\to
\Omega_F(\log)^{\oplus 2}$
sending $x$ to $(x,0)$
and to $(-x,x)$.
For both sides of
(\ref{eqmue}),
they are induced by the map
$\Omega_F(\log)
\to
\Omega_F(\log)^{\oplus 2}$
sending $x$ to $(-x,x)$.
\qed

{\it Proof of Lemma \ref{lmchi}.}
We may assume that
the residue field $F$
is separably closed.
Let $\chi$ be 
a character of
${\rm Gr}^rG_K$.
Let $\xi^*(\chi)$
be the character of
$G_{K'}^{\le r+}$
defined as the composition
of $\chi$ with
$$\begin{CD}
G_{K'}^{\le r+}
@>{\xi_*}>>
\pi_1(\Theta^{(r)}_{\bar F})^{\rm ab}
@>{(\ref{eqcan1})}>>
{\rm Gr}^rG_K.
\end{CD}$$
Since the canonical map
$\eta_*\colon
G_{K''}^{\le r+}
\to 
\pi_1(\Theta^{(r)}_{\bar F}
\times
\Theta^{(r)}_{\bar F})^{\rm ab}$
is surjective,
it suffices to show
the equality
${\rm p}_2^*\xi^*(\chi)
={\rm p}_1^*\xi^*(\chi)\cdot
\mu^*\xi^*(\chi)$
of characters of
$G_{K''}^{\le r+}$.

Since ${\rm Gr}^rG_K$
is a central subgroup of
$G^{\le r+}_K$,
there exists
an irreducible representation $V$
of $G^{\le r+}_K$
such that
the restriction to
${\rm Gr}^rG_K$
is the scalar
multiplication
by the character $\chi$.
By Corollary \ref{corchi},
we have an isomorphism
${\rm Res}_{1,\varepsilon}V
\to
{\rm Res}^{G_K^{\le r+}}
_{G_{K'}^{\le r+}}
V\otimes \xi^*(\chi)$
of representations
of $G_{K'}^{\le r+}$.
By Corollary \ref{cortra} and
by 
(\ref{eq2e}) and (\ref{eqmue}),
it induces
isomorphisms
\begin{align*}
{\rm Res}^{G_K^{\le r+}}
_{G_{K''}^{\le r+}}
V\otimes 
{\rm p}_2^*\xi^*(\chi)
&\to
{\rm Res}_{{\rm p}_1,\varepsilon'} 
{\rm Res}^{G_K^{\le r+}}
_{G_{K'}^{\le r+}}
V\otimes 
{\rm p}_1^*\xi^*(\chi),
\\
{\rm Res}^{G_K^{\le r+}}
_{G_{K''}^{\le r+}}
V \otimes 
\mu^*\xi^*(\chi)
&\to
{\rm Res}_{{\rm p}_1,\varepsilon'}
{\rm Res}^{G_K^{\le r+}}
_{G_{K'}^{\le r+}} V
\end{align*}
of representations
of $G_{K''}^{\le r+}$.
Thus, we obtain
an isomorphism
$${\rm Res}^{G_K^{\le r+}}
_{G_{K''}^{\le r+}}
V\otimes 
{\rm p}_2^*\xi^*(\chi)
\to
{\rm Res}^{G_K^{\le r+}}
_{G_{K''}^{\le r+}}
V \otimes 
{\rm p}_1^*\xi^*(\chi)
\cdot
\mu^*\xi^*(\chi).$$
Since
$G_{K''}^{\le r+}\to
G_K^{\le r+}$ is surjective,
this implies
an equality
${\rm p}_2^*\xi^*(\chi)
=
{\rm p}_1^*\xi^*(\chi)
\cdot
\mu^*\xi^*(\chi)$
of characters of
$G_{K''}^{\le r+}$
by Schur's lemma.
Thus the assertion is proved.
\qed

As in \cite[Corollary 1.26]{Jus},
Theorem \ref{them2}
has the following consequence.
Let $V$ be an
$\ell$-adic representation
of $G_K$.
Since $P=G_{K,\log}^{0+}$
is a pro-$p$ group,
there exists a unique
direct sum decomposition
$V=\bigoplus_{q\ge0,
q\in {\mathbb Q}}
V^{(q)}$
by sub $G_K$-modules
such that
the $G_{K,\log}^{r+}$-fixed part
is given by
$V^{G_{K,\log}^{r+}}=
\bigoplus_{q\ge r}
V^{(q)}.$
We put ${\rm Sw}_KV=
\sum_rr\cdot {\rm rank}\ V^{(r)}
\in {\mathbb Q}$.

\begin{cor}\label{corHA}
$${\rm Sw}_KV\in {\mathbb Z}[\frac1p].$$
\end{cor}

Liang Xiao claims 
a stronger assertion
${\rm Sw}_KV\in {\mathbb Z}$
in \cite[Theorem 3.5.11]{Xiao}.


\begin{thebibliography}{99}



\bibitem{AS1} A.\ Abbes and T.\ Saito, 
{\em Ramification of local fields
with imperfect residue fields I},
American J.\ of Mathematics,
124.5 (2002), 879-920.

\bibitem{AS2} ------,
{\em ibid.\ II},
Documenta Mathematica, Extra Volume Kato (2003), 3-70.


\bibitem{Jus} T.\ Saito, 
{\it Wild ramification and the characteristic cycle 
of an $\ell$-adic sheaf,}
Journal de l'Institut de Mathematiques de Jussieu, 
(2009) 8(4), 769-829 

\bibitem{Xiao} Liang Xiao,
{\it On Ramification Filtrations and p-adic Differential Equations, II: mixed characteristic case}, 
preprint, {\tt arXiv:0811.3792}

\end{thebibliography}
\end{document}